\documentclass[fleqn]{mat01}
\usepackage{times,mathtimy,amssymb,latexsym,amscd}
\begin{document}

\setcounter{page}{71}
\firstpage{71}

\newtheorem{defi}{\rm DEFINITION}
\newtheorem{propos}{\rm PROPOSITION}[ssection]
\newtheorem{rem}{Remark}
\newtheorem{corol}{\rm COROLLARY}[ssection]

\renewcommand{\theequation}{\arabic{section}.\arabic{equation}}

\def\d{\mbox{\rm d}}
\def\e{\mbox{\rm e}}


\title{Matrix multiplication operators on Banach function spaces}

\markboth{H Hudzik, Rajeev Kumar and Romesh Kumar}{Matrix multiplication operators on Banach function spaces}

\author{H HUDZIK, RAJEEV KUMAR$^{*}$ and ROMESH KUMAR$^{*}$}

\address{Faculty of Mathematics and Computer Science, Adam
Mickiewicz University Umultowska 87, 61--614, Poznan, Poland\\
\noindent $^{*}$Department of Mathematics, University of Jammu,
Jammu 180~006, India\\
\noindent E-mail: hudzik@amu.edu.pl; raj1k2@yahoo.co.in; romesh$_{-}$jammu@yahoo.com}

\volume{116}

\mon{February}

\parts{1}

\pubyear{2006}

\Date{MS received 23 August 2005; revised 27 September 2005}

\begin{abstract}
In this paper, we study the matrix multiplication operators on Banach
function spaces and discuss their applications in semigroups for solving
the abstract Cauchy problem.
\end{abstract}

\keyword{Banach function spaces;
closed operators; compact operators; Fredholm operators; matrix
multiplication operators; semigroups.}

\maketitle

\section{Introduction}

Let $(\Omega, \Sigma, \mu )$ be a $\sigma$-finite complete measure
space and $\mathbf{C}$ be the field of complex numbers.\@ By $L( \mu,
\mathbf{C}^{N} )$, we denote the linear space of all equivalence classes of
$\mathbf{C}^{N}$-valued $\Sigma$-measurable functions on $\Omega$
that are identified $\mu$-a.e. and are considered as column vectors.

Let $M_{\circ}$ denote the linear space of all functions in $L(\mu,
\mathbf{C}^{N} )$ that are finite a.e. With the topology of convergence
in measure on the sets of finite measure, it is a metrizable space.

The $\mathbf{C}^{N}$-valued Banach function space $X$ is defined as
\begin{equation*}
X= \{ f \in L( \mu, \mathbf{C}^{N} )\hbox{:}\ \| f \|_{X} < \infty \},
\end{equation*}
where $\|. \|_{X}$ is a function norm on $X$ such that for each $f, g,
f_n \in L( \mu, \mathbf{C}^{N} ), \, n \in \mathbf{N}$, we have

\begin{enumerate}
\leftskip .4pc
\renewcommand{\labelenumi}{(\roman{enumi})}
\item $0 \leq \| f (x) \|_{\mathbf{C}^N} \leq  \| g (x)
\|_{\mathbf{C}^N}$ for $\mu$-a.e. $x \in \Omega \Rightarrow \| f \|_{X} \leq \| g
\|_{X}$,

\item $0 \leq (f_n)_i  \nearrow (f)_i \; \mu$-a.e. for each
$i=1,2, \dots , N \Rightarrow \| f_n \|_{X}  \nearrow \| f \|_{X}$, and

\item $E \in \Sigma$ with $\mu (E) < \infty$ implies that
$ 1_{E}.z \in X$, for each $z \in \mathbf{C}^{N}$, and
\begin{equation*}
\hskip -1.25pc \int_E \|f (x) \|_{\mathbf{C}^{N} } \, \d \mu(x) \leq C_E \| f \|_{X},
\end{equation*}
for some constant $0 < C_E < \infty $, depending on $E$ but independent
of $f$, where $1_{E}$ is the characteristic function of the set $E$.
\end{enumerate}

\begin{defi}$\left.\right.$\vspace{.5pc}

\noindent {\rm A function $f$ in a Banach function space $X$ is said to have {\em
absolutely continuous norm} in $X$ if $ \| f 1_{E_{n} } \|_X \rightarrow 0
$ for every sequence  $ \{ E_n \}_{n=1}^{\infty }$ of $\mu$-measurable
sets in $\Omega$ satisfying $ E_n \rightarrow \emptyset , \;
\mu$-a.e., where $ E_n \rightarrow \emptyset$ means that $1_{E_{n}} \rightarrow
0, \; \mu$-a.e.}
\end{defi}

Let $X_a$ be the set of all functions in $X$ having {\it absolutely
continuous norm}. If $X_a =X$, then we say that $X$ has {\it absolutely
continuous norm}.

Let $X_b$ be the closure of the set of all $\mu$-simple functions in
$X$. Then, we have
\begin{equation*}
X_a \subseteq X_b \subseteq X.
\end{equation*}
Throughout this paper, we assume that $X= X_b$, that is, the simple
functions are dense in $X$. In case $X$ has {\em absolutely continuous
norm}, we have $X_a =X_b =X$ and so its Banach space dual $X^*$ and its
associate space $X'$ coincide, where $X'$ is defined as
\begin{equation*}
 X' = \{ g \in L(\mu, \mathbf{C}^N)'\hbox{:}\ \| g \|_{ X'} < \infty \}
\end{equation*}
and
\begin{equation*}
\| g \|_{X'} = \sup \left\lbrace \left\vert \int_{\Omega } \prec\! f(x),
g(x)\! \succ \, \d \mu (x) \right\vert\hbox{:}\  f \in X , \; \| f\|_X \leq 1 \right\rbrace,
\end{equation*}
where $L(\mu, \mathbf{C}^{N})'$ denotes the corresponding space of
equivalence classes of $\mathbf{C}^{N}$-valued functions considered as row
vectors. Note that $\prec\! f(x), g(x) \!\succ$ is the usual product of a row matrix
formed by $g(x)$ into a column matrix formed by $f(x)$.

The monotone convergence theorem holds in every Banach function space
$X$ in the form of the weak Fatou property (see axiom (ii)). Also note
that $X_a$ is the largest subspace of $X$ for which the suitable dominated
convergence theorem holds (see Proposition~3.6, p.~16 of \cite{bsi}).
So this fact can be used to easily generalise those results on
$L^p$-spaces to the general Banach function spaces having absolutely
continuous norm in which the dominated convergence theorem is required. It is
due to this fact that our results in \S3 follow on similar lines
as in $L^p$-spaces without any extra effort.

For details on Banach function spaces, we refer to \cite{bsi,linden,mali}.

For a measurable function $u\hbox{:}\  \Omega \rightarrow M_N(\mathbf{C})$, the
set of all $N \times N$ matrices over $\mathbf{C}$, the multiplication
transformation $M_u\hbox{:}\  L(\mu, \mathbf{C}^N) \mapsto L(\mu,
\mathbf{C}^N)$ is defined as
\begin{equation*}
M_u(f)= u \cdot f, \quad \hbox{for all}\  f \in L(\mu,
\mathbf{C}^N).
\end{equation*}

Using the arguments given in the introduction on p.~517 of \cite{siro} and p.~163, Theorem~1.2 of
\cite{enj} we can easily prove the next result.

\begin{theorem}[\!]
The multiplication operator $M_u$ is a bounded operator on a Banach
function space $X$ if and only if $u \in L^{\infty}(\mu, M_N(\mathbf{C}
))${\rm ,} the space of all $M_N(\mathbf{C})$-valued essentially-bounded
measurable functions. Moreover{\rm ,} we have
\begin{align}
 \| M_u \|_{X \mapsto X} & = \| u \|_{\infty} := \inf_{q \in [u] }
\sup_{x \in \Omega } \| q(x) \|,\\[.3pc]
\|M_{u}\|_{X\mapsto X} &= \inf \{M \geq 0\hbox{\rm :}\ \mu (\{ x \in \Omega\hbox{\rm :}\ \| u(x) \|
> M \} ) = 0 \},\hskip -1cm \phantom{0}
\end{align}
where $\| u(x) \|$ is the operator norm of $u(x)$ in $\mathbf{B}
(\mathbf{C}^N )$ induced by $\|\cdot\|_{\mathbf{C}^N }$ and $[u]$ is the
equivalence class of all measurable functions that are $\mu$-a.e. equal to $u$.
\end{theorem}

Note that for $N=1, \; L^{\infty}(\mu, M_N ( \mathbf{C} ))$ is the
space
of essentially bounded $\mathbf{C}$-valued measurable functions
denoted as $L^{\infty}(\mu)$.

For details on matrix analysis, see \cite{horn}.

\begin{defi}$\left.\right.$\vspace{.5pc}

\noindent {\rm For $u \in L  (\mu, M_N (\mathbf{C}))${\rm ,}
an operator $(M_u, D(M_u))$ defined on $X( \mathbf{C}^N)$ by
$(M_u f(x)) = u(x) \cdot f(x)$ for $x \in \Omega$ and for all $f \in D(M_u)= \{ f
\in X\hbox{:}\ u\cdot f \in X \}$ is called a matrix multiplication operator.}
\end{defi}

We consider an abstract Banach space-valued linear initial value
problem of the form
\begin{align*}
  \dot{v}(t) & = A v(t),\quad t \geq 0, \\[.3pc]
   v(0)      & =  x,
 \end{align*}
where the independent variable $t$ represents time, $v(\cdot)$ is a
function with values in a $\mathbf{C}^{N}$-valued Banach function space
$X=X(\mathbf{C}^{N } )$, where $\mathbf{C}^{N}$ denotes the
$N$-dimensional complex space. Since all the norms on $\mathbf{C}^{N}$
are equivalent, we choose $\|\cdot\|_{\mathbf{C}^{N}}=
\|\cdot\|_{\sup}$. Note that this norm induces the matrix norm
$\|(a_{ij})_{N \times N} \|= \max_{i=1, \dots, N} \sum_{j=1}^{N} \vert
a_{ij} \vert$ on $M_N (\mathbf{C})$, that is, the norm of a matrix is
given by the maximal absolute sum of its rows.

Also, $A\hbox{:}\ D(A) \subseteq X \mapsto X$ is a linear operator and $x \in
X$ is the initial value. This problem is called an abstract Cauchy
problem (ACP) associated with a matrix multiplication operator $(M_u,
D(M_u))$ and the initial value $x$.

A function $u\hbox{:}\ \mathbf{R}_{+} \mapsto X$ forms a solution of ACP, if $u$
is a continuously differentiable function with
respect to $X$ and $ u(t)\in D(A)$, for each $t \geq 0$ (see
\cite{nue}).

If the operator $A$ is the generator of a strongly continuous semigroup
$(T(t))_{t \geq 0}$, then for each $x \in D(A)$, the function
\begin{equation*}
u\hbox{:}\ t \mapsto u(t)=T(t)x
\end{equation*}
is the unique solution of ACP. See \cite{nue} and \cite{pazy} for
details on ACP and semigroups. Note that every well-posed Cauchy
problem is solved by a strongly continuous semigroup and each such
semigroup is the solution semigroup of a well-posed Cauchy problem (see
\cite{nag2}).

In this paper, we extend the results of \cite{enj,holder,komal,takagi}
to general Banach function spaces. Note that the separable Banach
function spaces form a subclass of the absolutely continuous ones. The
examples of Banach function spaces having absolutely continuous norm are
$L^p$-spaces, Orlicz spaces with $\Delta_2$-conditions \cite{rao},
Lorentz spaces \cite{bsi}, separable Orlicz--Lorentz spaces
\cite{hudzik}, etc. In the second section, we study the compactness,
closedness, invertibility and Fredholmness properties of a
multiplication operator $M_u$ on a Banach function space $X$. In the
third section, we discuss the abstract Cauchy problem associated with
the operator $(M_u, D(M_u))$ on $X$.

\section{Properties of multiplication operators}

\setcounter{equation}{0}

For a measurable function $u\hbox{:}\ \Omega \mapsto M_N(\mathbf{C})$,
we call the set
\begin{equation*}
u_{\rm ess}(\Omega) = \{ \lambda \in \mathbf{C}\hbox{:}\ \mu ( \{s \in
\Omega\hbox{:}\ \| u(s)- \lambda \| < \epsilon \} ) \neq 0,\quad \forall
\epsilon >0 \},
\end{equation*}
its essential range and define the associated multiplication operator
$M_u$ on the space $X$ by $M_u(f)=u\cdot f$ for each $f$ in the domain
\begin{equation*}
D(M_u)= \{ f \in X\hbox{:}\ u\cdot f \in X \}.
\end{equation*}
In case $M_u$ is a bounded linear operator on $X$, we have $D(M_u)=X$.

\begin{propos}\label{th:one}$\left.\right.$\vspace{.5pc}

\noindent Let $(M_u, D(M_u))$ be the multiplication operator on a Banach
function space $X$ induced by some measurable function $u\hbox{\rm :}\ \Omega
\mapsto M_N (\mathbf{C})$. Then the following statements hold{\rm :}

\begin{enumerate}
\leftskip .4pc
\renewcommand{\labelenumi}{\rm (\roman{enumi})}
\item The operator $(M_u,D(M_u))$ is closed and densely defined{\rm ,}
in case $X$ has absolutely continuous norm.

\item The operator $M_u$ has a bounded inverse if and only if
$0 \notin u_{ {\rm ess} } (\Omega)$.
In this case{\rm ,} $M_{u^{-1}}= M_r${\rm ,} for some measurable
$r\hbox{\rm :}\ \Omega \mapsto M_N (\mathbf{C})$ defined by
\begin{equation*}
\hskip -1.25pc r(s)= \begin{cases}
(u(s))^{-1}  & if\ u(s)\neq 0,\\
0            & if\ u(s)=0,\end{cases}
\end{equation*}
for each $s \in \Omega$.

\item The spectrum of $M_u$ is the essential range of $u${\rm ,}
i.e.{\rm ,} $\sigma (M_u) = u_{\rm ess} (\Omega )$.
\end{enumerate}
\end{propos}

The proof follows on similar lines as in $L^p$-spaces (see \cite{enj}
and \cite{nagel}).

\begin{theorem}[\!]\label{th:algebra}
For $N=1${\rm ,} the set of all bounded multiplication operators on a Banach
function space $X$ of~\,$\mathbf{C}$-valued functions forms a maximal
abelian subalgebra of~\,$\mathbf{B}(X)${\rm ,} the space of all bounded linear
operators on $X$.
\end{theorem}

\begin{proof}
Let $\mathbf{m} =\{ M_u\hbox{:}\ u \in L^{\infty}(\mu ) \} $. Clearly,
$\mathbf{m}$ is an abelian subalgebra of $\mathbf{B}(X)$. We prove that
$\mathbf{m}$ is maximal, i.e., if $A$ commutes with $\mathbf{m}$, then
$A \in \mathbf{m}$.

Let $e\hbox{:}\ \Omega \to \mathbf{C}$ be the unity function. Let $v=Ae$
and $E \in \Sigma$. Then
\begin{equation*}
A 1_E = A M_{ 1_E } e = M_{ 1_E } A e = 1_E   v = v 1_E =     M_v
1_E .
\end{equation*}
We claim $v\in L^{\infty}(\mu )$. Suppose that the set
\begin{equation*}
F_n =\{ x \in\Omega\hbox{:}\ |v(x)| > n \}
\end{equation*}
has a positive measure for each $n \geq 1$. By the finite subset property
of the underlying measure space, we assume $\mu ( F_n )< \infty$. Then,
$1_{F_n} \in X$ and we have
\begin{align*}
n 1_{F_n} (x)
& \leq 1_{F_n} |v(x)|\\[.2pc]
\Rightarrow n \| 1_{F_n} \|_{X}
& \leq \| M_v 1_{F_n} \|_{X}  \\[.2pc]
\Rightarrow n \| 1_{F_n} \|_{X}
& \leq \| A 1_{F_n} \|_{X}
\end{align*}
for each $n \in \mathbf{N}$, which contradicts the boundedness of $A$.
Therefore $v \in L^{\infty} (\mu )$. Since the
set of $\mu$-simple functions is dense in $X$, we have
$A= M_v$. This proves that $A \in \mathbf{m}$ and so $\mathbf{m}$ is
a maximal abelian subalgebra of $ \mathbf{B}(X)$.
\end{proof}

\begin{theorem}[\!]\label{th:closed}
Let $X=X(\mathbf{C})$ be a Banach function space of
$\mathbf{C}$-valued measurable functions on $\Omega$ and $M_u \in \mathbf{B} (X)$ for some
$u \in L^{ \infty}(\mu)$. Then $M_u$ has closed
range if and only if there exists some $ \delta >0 $ such that
$|u(x)| \geq \delta, \; {\rm for} \; \mu $-almost all $x \in
{\rm support}(u)=S$.
\end{theorem}

\begin{proof}
Suppose $|u(x)|\geq \delta $ for $\mu$-almost all $x \in S$.
Then using the same technique as in converse part of Theorem~$1.2$ in
p.~163 of \cite{enj} with $N = 1$, we have
\begin{equation*}
|(M_u\big|_S f)(x)|\geq \delta |f(x)|,
\end{equation*}
for $\mu$-a.e. $x \in \Omega$ and each $f\in X$. As the norm
$\|\cdot \|_{X}$ is increasing on $X$, so we have
\begin{equation*}
\| M_u |_S f\|_{X} \geq \delta \| f \|_{X},
\end{equation*}
for each $f\in X$. This implies that $M_u$ is invertible. Thus, we
conclude that $M_u$ has closed range.

Conversely, suppose $M_u$ has closed range. Then there exists some $
\delta >0$ such that
\begin{equation}\label{eq:delta}
\| M_u h \|_{X} \geq \delta \| h \|_{X}, \; {\rm for} \; {\rm each} \;
h \in X.
\end{equation}
Let $E=\{ x\in \Omega\hbox{:}\ |u(x)| < \delta / 2 \} $ be
such that $ \mu(E) >0$. So, there exists a measurable set $F \subseteq
E$ such that $1_F \in X $. Further,
\begin{equation*}
|(M_u 1_F)(x)| = \vert u(x) 1_F (x) \vert <\delta
|1_F(x)|,
\end{equation*}
implies that
\begin{equation*}
\| M_u 1_F \|_X \leq \delta \| 1_F \|_X,
\end{equation*}
which contradicts (\ref{eq:delta}) and so $\mu (E)=0$. Therefore $\vert
u(x)\vert \geq \delta / 2$, for $ \mu $-almost all  $x \in S$.
\end{proof}

\begin{theorem}[\!]
Let $X =X (\mathbf{C}^N)$ be a Banach function space of $\mathbf{C}^N$-valued measurable functions on $\Omega$ and $M_u \in  \mathbf{B}(X)$. Then
$M_u$ is a compact operator if and only if $X( N, \epsilon ,
\mathbf{C}^N)$ is finite-dimensional{\rm ,} for each  $\epsilon > 0 ${\rm ,} where
\begin{equation*}
N = N (u, \epsilon) = \{x \in \Omega\hbox{\rm :}\ \| u(x) \| \geq \epsilon \}
\end{equation*}
and
\begin{equation*}
X (N ,\epsilon , \mathbf{C}^N)= \{ f \in X\hbox{\rm :}\  f(x)= \mathbf{0}  \;{\rm
for}\; x \notin N \}.
\end{equation*}
\end{theorem}

\begin{proof}
Suppose $M_u$ is a compact operator. Then, its
restriction to the invariant subspace $X(N, \epsilon , \mathbf{C}^N)$ is
compact. Also, $M_u \big|_{ X(N, \epsilon , \mathbf{C}^N)}$ has closed range in
$X(N, \epsilon, \mathbf{C}^N)$. But $M_u \big|_{X( N, \epsilon ,
\mathbf{C}^N)}$ is invertible. Therefore,
$X( N, \epsilon , \mathbf{C}^N)$ is finite-dimensional for each
$\epsilon >0$.

Conversely, suppose $X( N, \epsilon , \mathbf{C}^N)$ is
finite-dimensional for each  $\epsilon >0$.
In particular, $X (N ,1/n , \mathbf{C}^N)$ is finite-dimensional for
each $n \geq 1$.  Take $w \in L^\infty( \mu, M_N( \mathbf{C}))$ and define $w_n\hbox{:}\ \Omega
\to M_N( \mathbf{C})$, by
\begin{equation*}
w_{n,ij} (x)=  \begin{cases}
w_{ij}(x)   & \mbox{ if}\ x \in A_n ,\\[.2pc]
 0       & \mbox{ if}\ x \notin A_n ,\end{cases}
\end{equation*}
for each $i, j =1, 2, \dots, N$,
where $ A_n = \{ x \in \Omega\hbox{:}\ \| w(x) \| \geq 1/n \} $.

Then for each $f \in X$, we have
\begin{align*}
{\| ((M_{w_n} - M_w)f)(x) \|_{\mathbf{C}^N} }  &\leq \| (w_n(x) - w(x) )f(x) 1_{A_n} (x) \|_{\mathbf{C}^N}\\[.2pc]
&\quad\,  + \| (w_n(x) -w(x))f(x) 1_{\Omega \backslash A_n} (x)\|_{
\mathbf{C}^N}\\[.2pc]
&=  \| w(x) f(x) 1_{\Omega \backslash A_n} (x)\|_{ \mathbf{C}^N} \leq  \frac{1}{n} \| f(x) \|_{ \mathbf{C}^N},
\end{align*}
for $\mu$-a.e., $x \in \Omega$. Thus, for each $f \in X$, we have
\begin{equation*}
\| (M_{w_n} -M_w)f\|_{X} \leq 1/n \| f \|_{X}.
\end{equation*}

This implies that $M_{w_n} \to M_w$ as $n \to \infty $. But each
$M_{w_n} $ is of finite rank and thus $M_w$ is a compact operator.
\end{proof}

The next result gives the necessary and sufficient conditions for a
multiplication operator $M_u$ on a Banach function space
$X=X(\mathbf{C})$ to be a Fredholm operator thereby generalising the results in
\cite{komal} for Orlicz spaces and \cite{takagi} for $L^p$-spaces. Here, we take $N=1$.

\begin{theorem}[\!]
Suppose $( \Omega , \Sigma, \mu)$ is a non-atomic measure space and
$M_u \in \mathbf{B}(X)${\rm ,} where $X=X(\mathbf{C} )$ is a Banach
function space having absolutely continuous norm. Then{\rm ,} the following are
equivalent.

\begin{enumerate}
\leftskip .4pc
\renewcommand{\labelenumi}{\rm (\roman{enumi})}
\item $M_u$ is invertible.
\item $M_u$ is Fredholm.
\item The range of $M_u${\rm ,} that is{\rm ,} $R(M_u)$ is closed and
${\rm codim}(R(M_u))< \infty$.
\item $|u(x)|\geq \delta $ for $\mu$-a.e. $x \in \Omega$
for some $ \delta > 0 $.
\end{enumerate}
\end{theorem}

\begin{proof}
We prove (iii) $\Rightarrow$ (iv) only, as the other implications are
obvious by using the previous results. Suppose $R(M_u)$ is closed and
${\rm codim}(R(M_u))< \infty$.

We claim $M_u$ is onto. Suppose the contrary. Then there exists
$f_{\circ} \in X\backslash R(M_u)$. Since $R(M_u)$ is closed, there exists
$g_{\circ} \in X^*, $ the dual (associate) space of $X$ such that
\begin{equation}\label{eq:integer}
\int_{\Omega } \prec \!f_{\circ}, g_{\circ} \!\succ \, \d \mu = 1
\end{equation}
and
\begin{equation}\label{eq:inte}
\int_{\Omega } \prec \!(M_u f ), g_{\circ} \!\succ \, \d \mu = 0 , \; {\rm
for}\; {\rm each} \; f\in X.
\end{equation}

Now (\ref{eq:integer}) yields that the set
\begin{equation*}
E_{\delta}=\{ x \in \Omega\hbox{:}\ \hbox{Re} \, (  \prec\! f_{\circ}(x),
g_{\circ}(x) \!\succ) \geq \delta \}
\end{equation*}
has positive $\mu$-measure for some $\delta >0$. As $\mu $ is
non-atomic, we can choose a sequence $\{ E_n \}$ of subsets of $E_{\delta}$
with $ 0< \mu (E_n) < \infty  \; {\rm and}\; E_m \cap E_n = \emptyset$ for
$m \neq n$.

Put $g_n= 1_{E_n} g_{\circ}$.
Clearly, $g_n \in X^*, \;{\rm as}\; \|  g_n \|_{X^*} \leq \| g_{\circ}
\|_{X^*} \; {\rm and} \; g_n \neq 0$, as
\begin{equation*}
{\rm Re} \, \int_{\Omega} \prec\! f_{\circ}, g_n \!\succ \, \d \mu
= {\rm Re} \int_{E_n} \prec\! f_{\circ}, g_{\circ} \!\succ \, \d \mu \geq
\delta \mu (E_n) > 0 ,
\end{equation*}
for each $n$. Also, for each $f  \in X , \; 1_{E_n } f \in X$ and
so (\ref{eq:inte}) implies that
\begin{align*}
\int_{\Omega } \prec\! f, ( M_u^* g_n) \!\succ  \, \d \mu
& = \int_{\Omega } \prec\! M_u f,g_n \!\succ \, \d \mu
= \int_{\Omega } \prec\! M_u f, 1_{E_n} g_{\circ} \!\succ \, \d \mu  \\[.3pc]
& =  \int_{\Omega } \prec\! M_u 1_{E_n} f,  g_{\circ} \!\succ \, \d \mu
=0,
\end{align*}
where $M_u^*$ is the conjugate operator of $M_u $, which implies that
$M_u^* g_n =0 $ a.e. and so $g_n\in {\ker}(M_u^*)$. Since all the sets
in  $\{E_n\}$ are disjoint, the sequence $\{g_n\}$ forms a
linearly-independent subset of $\ker(M_u^*)$. This contradicts the
fact that
\begin{equation*}
\dim \, \ker (M_u^*)= \hbox{codim} \, R(M_u) < \infty.
\end{equation*}
So $M_u$ is onto.

Let $\ker(u)= \{ x \in \Omega\hbox{:}\ u(x) = 0 \}$.
Then $\mu ({\ker}(u))=0$. Since $\mu (\ker(u))>0$, there is an
 $ F \subseteq \ker (u) \; {\rm with}\; 0< \mu (F)<\infty$.
Thus, $1_F \in X\backslash R(M_u)$ which contradicts the
surjectiveness of $M_u$. For $n\geq 1$, put
\begin{equation*}
G_n = \left\lbrace x \in \Omega\hbox{:}\ \frac{ \| u \|_{\infty} }{ (n+1)^2} < |u(x)|
\leq \frac{ \| u \|_{\infty} }{ n^2}\right\rbrace
\end{equation*}
and
\begin{equation*}
 T =\{ n\in \mathbf{N}\hbox{:}\ \mu (G_n) > 0 \}.
\end{equation*}
Clearly $ \Omega = \cup_{n=1}^{\infty} G_n $ and $ \mu(G_n) <\infty$,
for each $n\geq 1$. Take
\begin{equation*}
f = \sum_{ n \in T } \frac { u \cdot 1_{G_{n} }  }{ \| 1_{ G_{n} } \|_{X}} .
\end{equation*}
Then
\begin{align*}
|f(x) | & =  \left\vert \sum_{n \in T} \frac{ u(x) 1_{G_n}(x)  }{ \| 1_{G_n}
\|_{X}  } \right\vert\\[.3pc]
& \leq \Sigma_{n\in T} \frac{ | u(x)  1_{G_n}(x)  |
 }{ \| 1_{G_n} \|_{X} } \\[.3pc]
&\leq \Sigma_{ n \in T} \frac{  \| u \|_{ \infty }  }{ n^2  \| 1_{G_n} \cdot z
\|_{X} } ,
\end{align*}
for each $x \in \Omega $, this proves that
\begin{equation*}
\| f \|_{X} \leq \sum_{n \in T} \frac{ \|u \|_{\infty }  }{n^2  }
\leq \|u \|_{\infty } \sum_{n=1 }^{\infty} \frac{1}{n^2}
< \infty.
\end{equation*}
Therefore $f\in X$ and so there exists some $g \in X$ such that $M_u \,
g=f$
with
\begin{equation*}
 g = \frac{f }{ u} = \Sigma_{n\in T} \frac{1_{G_n}  }{ \|   1_{G_n} \|_{X} } .
\end{equation*}
Since the sets $G_n$ are disjoint, using the absolute continuity of $X$
(Theorem~4.1, p.~20 of \cite{bsi}), we have
\begin{align*}
\| g \|_{X}
&  = \sup_{h \in X^*, \|h \|_{X^*} \leq 1}\left\vert \int_{\Omega} \prec\!
\sum_{n \in T} \frac{ 1_{G_n}  }{ \|   1_{G_n} \|_{X} },  h \succ \,
\d \mu \right\vert\\[.3pc]
&  = \left( \sum_{n \in T} \frac{1 }{ \|  1_{G_n} \|_{X} } \right) \sup_{h \in
X^*, \|h \|_{X^*} \leq 1}\left\vert \int_{\Omega} \prec\! 1_{G_n} , h \!\succ
\, \d \mu \right\vert\\[.3pc]
&  = \sum_{n \in T} 1 .
\end{align*}
Then, $\| g \|_{X} < \infty$ if and only if $T$ is finite. So
there is some $m>0$ such that for $n \geq m , \; \mu(G_n)=0 $ and
$ \mu({\ker}(u))=0$, which further implies that
\begin{align*}
\mu \{ x \in \Omega\hbox{:}\ \| u(x) \| \leq \|u\|_{\infty}/m^2 \}
&= \mu((\cup_{n=m}^{\infty} G_n) \cup \ker(u))  \\[.3pc]
&\leq \mu((\cup_{n=1}^{\infty} G_n) \cup \ker(u)) = 0,
\end{align*}
that is,
$ \| u(x) \| \geq  \| u \|_{\infty}/ m^2 = \delta $, a.e. on $X$. This
proves (iv).
\end{proof}

\section{Semigroups of multiplication operators}

\setcounter{equation}{0}

In this section, we study some applications of the multiplication
operators on $X$ in semigroup theory. Using Proposition~\ref{th:one}(i),
we see that $M_u$ is a closed and densely defined multiplication operator
defined on an absolutely continuous Banach function space $X$.

The essential spectrum of a multiplication operator is defined as
\begin{equation*}
\sigma_{\rm ess} (M_u)= \cup_{x \in \Omega } \sigma_{\rm ess} (u(x))= \cap_{ p \in [u] }
\overline{ \cup_{ x \in \Omega } \sigma (p(x)) },
\end{equation*}
where $\sigma (p(x))$ denotes the spectrum of the matrix $p(x)$ and $[u]$ is the
equivalence class of all measurable functions that are $\mu$-a.e. equal to $u$.

For an open $\epsilon$-disk $U_{\epsilon}$
with center $0$, we have
\begin{align*}
\sigma_{\rm ess} (M_u)&= \cup_{x \in X} \sigma_{\rm ess} (u(x))\\[.3pc]
&=\{ z \in \mathbf{C}\hbox{:}\ \mu( \{
x \in \Omega\hbox{:}\ \sigma (u(x)) \cap z+U_{\epsilon} \} )>0,\quad \forall
\epsilon > 0   \}.
\end{align*}
Using Propositions~4.11 and 4.12, p.~32 of \cite{nagel}, it
is easy to prove the result.

\setcounter{propos}{0}
\begin{propos}$\left.\right.$\vspace{.5pc}

\noindent Suppose $(M_u , D(M_u))$ is a matrix multiplication operator on the
Banach function space $X$ with non-void resolvent set $\rho (M_u)$. Then
its spectrum is given by
\begin{equation*}
\sigma(M_u)=\overline{\cup_{x \in \Omega} \sigma_{\rm ess} (u(x))} .
\end{equation*}
\end{propos}

\begin{rem}
{\rm If $\Omega= \mathbf{R}^{m}$ with the Lebesgue measure $\mu$ and $u$ is
a continuous function with non-void resolvent, then we have
\begin{equation*}
\sigma(M_u)= \overline{ \cup_{x \in \Omega} \sigma (u(x))}.
\end{equation*}

The stability of the solutions of the abstract Cauchy problem are
determined by the spectral bound of the corresponding operator, \cite{enj}.}
\end{rem}

\setcounter{corol}{1}
\begin{corol}$\left.\right.$\vspace{.5pc}

\noindent A matrix multiplication operator $M_u$ is bounded if and only if its
spectrum $\sigma (M_u)$ is bounded.
\end{corol}

There are a number of functional analytic approaches to the
well-posedness of solutions of abstract Cauchy problems (see, for example,
\cite{bb,nagel,gold,holder}).
The proofs of the following results are on the similar lines as in
\cite{holder} for $L^p$-spaces, so we only state the results.

\setcounter{defin}{2}
\begin{theorem}[\!]
A matrix multiplication operator $M_u$ generates a semigroup on $X$
having absolutely continuous norm if and only if
\begin{equation}
\overline{ \lim_{t \to 0}  } \|  \e^{tu} \|_{\infty} < \infty.
\end{equation}

In case $M_u$ generates the semigroup $\{ T(t) \}_{t \geq 0}${\rm ,} then we
have $T(t)=M_{{\rm e}^{tu}}$. On the other hand{\rm ,} if $\{ M_{{\rm e}^{tu}}
\}_{t \geq 0}$ defines a semigroup on $X${\rm ,} then its generator is
given by $M_u $.
\end{theorem}

\begin{theorem}[\!]
If the matrix multiplication operator $(M_u, D(M_u))$ generates the
semigroup $\{ T(t) \}_{t \geq 0}${\rm ,}
then
\begin{equation*}
\sigma(T(t))= \overline{ \e^{t \sigma(M_u)}}.
\end{equation*}
\end{theorem}

\begin{defi}{\rm \cite{holder,linden}}$\left.\right.$\vspace{.5pc}

\noindent {\rm Let $(A,D(A))$ be a linear operator on a Banach space $F$
such that there exists some constants  $ m \in \mathbf {N} , M >0 , \,
w \in \mathbf{R}$, the real line and a strongly continuous family
$(S(t))_{t \geq 0}$ in $\mathbf{B} (F)$, the space of all the bounded linear
operators on $F$ with
\begin{equation*}
\| S(t) \| \leq M\e^{wt}, \; {\rm for} \; {\rm each} \; t \geq 0
\end{equation*}
such that its resolvent operator $R(\lambda, A)= (\lambda - A)^{-1}$
exists and is given by
\begin{equation*}
 R(\lambda, A) y= \lambda^{m} \int_{0}^{ \infty } \e^{- \lambda t} S(t)
y \, \d t ,
\end{equation*}
for $y \in F$ and $\lambda > w$. Then $(A, D(A))$ is the generator of
an $m$-times integrated semigroup $(S(t))_{t \geq 0}$.

The generator of an $m$-times integrated semigroups yields the
well-posedness of the abstract Cauchy problem, since there exists some
constants $M$ and $w$ as above such that for all initial values $u_{\circ}
\in D(A^{m+1})$, we get a unique solution $u(\cdot)$ with the\break property
\begin{equation*}
\| u(t) \| \leq M\e^{wt} \|x\|_{A^m},
\end{equation*}
for each $t \geq 0$, where $\|\cdot\|_{A^m}$ denotes the $m$th-graph
norm of $(A, D(A))$ (see \cite{nue}).}
\end{defi}

\begin{theorem}[\!]
Let $(M_u, D(M_u))$ be a matrix multiplication operator on a Banach
function space $X= X(\mathbf{C}^{N})$. Then the following are equivalent.

\begin{enumerate}
\leftskip .4pc
\renewcommand{\labelenumi}{\rm (\roman{enumi})}
\item The operator $(M_u, D(M_u))$ is the generator of an
integrated semigroup.
\item There exists a constant $w \in \mathbf{R}$ such that
\begin{equation*}
\hskip -1.25pc \sigma (M_u) \subseteq \{ z \in \mathbf{C}\hbox{:}\ {\rm Re} \ z \leq w  \}.
\end{equation*}
\item The resolvent set of $M_u$ is non-empty such that
\begin{equation*}
\hskip -1.25pc {\rm ess} \sup_{x \in \Omega} s(u(x)) \leq w,
\end{equation*}
for some $w \in \mathbf{R}${\rm ,} where $s(u(x))$ denotes the spectral bound
of the matrix $u(x) \in M_N (\mathbf{C} )$.
\end{enumerate}
\end{theorem}

\begin{rem}
{\rm As in Corollary~4.9 of \cite{nue}, we see that the multiplication
operator $(M_u, D(M_u))$ satisfying one of the conditions in the statement of
the above-stated theorem is a generator of a $(2N+1)$-times integrated
semigroup.

For a multiplication operator $(M_u, D(M_u))$ on $X=X(\mathbf{C}^{N})$,
we see that it is a generator of a strongly continuous semigroup if and
only if there is some  $c>0$ such that
\begin{equation*}
   \| \e^{tu } \|_{\infty} <c, \;  {\rm for}\; t \in [0, 1].
\end{equation*}}
\end{rem}

\begin{rem}
{\rm Theorem~1 on p.~16 of \cite{holder} is also true for Banach
function spaces $X=X(\mathbf{C}^{N})$. Further, if $X$ has absolutely
continuous norm, then p.~163, Theorem~2 of \cite{holder}, p.~164,
Proposition~1 of \cite{holder} and p.~165, Proposition~2 of \cite{holder}
can be proved for such Banach function spaces.}
\end{rem}

\begin{rem}
{\rm The matrix multiplication operator $(M_u, D(M_u))$ on $X$ also
generates some subclasses of the strongly continuous semigroups on $X$
such as analytic, differentiable, norm-continuous semigroups, etc. See
\cite{nagel,gold,pazy} for details on these types of semigroups.\pagebreak

\noindent Also,
we see that the analyticity of semigroups on $X$ depends on the spectrum
$\sigma (M_u)$ of the matrix multiplication operator $(M_u, D(M_u))$ on
$X$.}
\end{rem}

\section*{Acknowledgements}

The authors are thankful to the referee for valuable suggestions and
comments. One of the authors (Rajeev Kumar) is supported by CSIR-grant
(sanction no. 9(96)/2002-EMR-I, dated 13-5-2002).

\end{document}